\documentclass[11pt, final]{article}
 \usepackage{a4}
 \usepackage{amsmath}%
 \usepackage{amstext}%
 \usepackage{amssymb}%
 \usepackage{showkeys}%
\usepackage{cite}
\newcommand{\R}{\mathbb{R}}%
\DeclareMathOperator*\inte{int}%
\DeclareMathOperator*\PMin{PMin}%
\DeclareMathOperator*\Min{Min}%
\DeclareMathOperator*\Max{Max}%

\title{Classical linear vector optimization duality revisited
\thanks{Research partially supported by DFG (German Research Foundation), project WA 922/1-3.}}

\author{Radu Ioan Bo\c t
\thanks {Faculty of Mathematics, Chemnitz University of Technology,
D-09107 Chemnitz, Germany, e-mail:
radu.bot@mathematik.tu-chemnitz.de.}, Sorin-Mihai Grad \thanks
{Faculty of Mathematics, Chemnitz University of Technology,
D-09107 Chemnitz, Germany, e-mail:
sorin-mihai.grad@mathematik.tu-chemnitz.de.}, and
Gert Wanka \thanks
{Faculty of Mathematics, Chemnitz University of Technology,
D-09107 Chemnitz, Germany, e-mail:
gert.wanka@mathematik.tu-chemnitz.de.}
}
\date{}
\begin{document}
\maketitle
\textbf{Abstract.} With this note we bring again into attention a vector dual problem neglected by the contributions who have recently announced the successful healing of the \textit{trouble} encountered by the classical duals to the classical linear vector optimization problem. This vector dual problem has, different to the mentioned works which are of set-valued nature, a vector objective function.
Weak, strong and converse duality for this ``new-old'' vector dual problem are proven and we also investigate its connections to other vector duals considered in the same framework in the literature. We also show that the efficient solutions of the classical linear vector optimization problem coincide with its properly efficient solutions (in any sense) when the image space is partially ordered by a nontrivial pointed closed convex cone, too.\\


\section{Introduction and preliminaries}

All the vectors we use in this note are column vectors, an upper index ``$^T$'' being used to transpose them into row vectors. Having a set $S\subseteq \R^k$, by $\inte(S)$ we denote its interior. By $\R^k_+$ we denote the \textit{nonnegative orthant} in $\R^k$. We say that a set $K\subseteq \R^k$ is a \textit{cone} if $\lambda K\subseteq K$ for all $\lambda \in \R_+$. Then $K$ induces on $\R^k$ a partial ordering ``$\leqq_K$'' defined by $v\leqq _K w$ if $w-v\in K$. If $v\leqq_K w$ and $v\neq w$ we write $v\leq _K w$. When $K=\R^k_+$ these cone inequality notations are simplified to ``$\leqq$'' and, respectively, ``$\leq$''. A cone $K\subseteq \R^k$ which does not coincide with $\{0\}$ or $\R^k$ is said to be \textit{nontrivial}. A cone $K$ is called \textit{pointed} if $K\cap (-K) = \{0\}$. The set $K^*=\{\lambda\in \R^k: \lambda^Tv\geq 0 \ \forall v\in K\}$ is the \textit{dual cone} of the cone $K$. The \textit{quasi-interior} of $K^*$ is the set $K^{*0}= \{\lambda\in \R^k: \lambda^Tv> 0 \ \forall v\in K\backslash \{0\}\}$. The \textit{recession cone} of a convex set $M\subseteq \R^k$ is $0^+M=\{x\in \R^k: M + x \subseteq M\}$. A set is said to be \textit{polyhedral} if
it can be expressed as the intersection of some finite collection of closed half-spaces. For other notions and notations used in this paper we refer to \cite{RC}.

The vector optimization problems we consider in this note consist of vector-minimizing or vector-maximizing a vector function
with respect to the partial ordering induced in the image space of the vector function by a nontrivial pointed closed convex cone. For the
vector-minimization problems we use the notation $\Min$, while the vector-maximization ones begin with $\Max$. The solution concepts considered for these problems are based on the following minimality concepts for sets.

Let the space $\R^k$ be partially ordered by a nontrivial pointed closed convex cone $K\subseteq \R^k$ and $M\subseteq \R^k$ be a nonempty set.
An element $\bar x \in M$ is said to be a \textit{minimal element of $M$ (regarding the partial ordering
induced by $K$)} if there exits no $v \in M$ satisfying $v \leq_K \bar v$. The set of all minimal elements of $M$ is denoted by $\Min(M, K)$.
Even if in the literature there are several concepts of proper minimality for a given set (see \cite[Section 2.4]{DVO} for a review on this subject), we deal here only with the properly minimal elements of a set in the sense of linear scalarization. Actually, all these proper minimality concepts coincide when applied to the primal classical linear vector optimization problem we treat in this note, as we shall see later. An element $\bar v \in M$ is said to be a \textit{properly minimal element of $M$ (in the sense of linear scalarization)} if there exists a $\lambda\in K^{*0}$ such that $\lambda^T\bar v \leq \lambda ^Tv$ for all $v \in M$. The set of all properly minimal elements of $M$ (in the sense of linear scalarization) is denoted by $\PMin(M, K)$. It can be shown that every properly minimal element of $M$ is also minimal, but the reverse assertion fails in general. Corresponding maximality notions are defined by using the definitions from above. The elements of the set $\Max (M, K):= \Min (M, -K)$ are called \textit{maximal elements of $M$}.\\

The \textit{classical linear vector optimization problem} is
$$\begin{array}{ll}
(P)\hspace{0.2cm} & \Min\limits_{x \in {\cal A}} Lx,\\
& {\cal A}= \{x \in \R^n_+ : Ax=b\}
\end{array}$$
where $L\in \R^{k\times n}$, $A\in \R^{m\times n}$, $b\in \R^m$ and the space $\R^k$ is partially ordered by the nontrivial pointed closed convex cone $K\subseteq \R^k$.
The first relevant contributions to the study of duality for $(P)$ were brought by Isermann in \cite{ISZ, IOR, IZOR} for the case $K=\R^k_+$. The dual he assigned to it,
$$(D^I)\hspace{0.2cm}\Max\limits_{U \in {\cal B}^I} h^I(U),$$
where
$${\cal B}^I = \big \{U \in \R^{k\times m}
: \nexists x\in \R^n_+ \mbox{ such that } (L-UA)x\leq 0\big \}$$ and
$$h^I(U) = Ub,$$
turned out to work well only when $b\neq 0$ (see \cite{JMP}).

The same drawback was noticed in \cite{JMP, JB} also for the so-called \textit{dual abstract optimization problem} to $(P)$
$$(D^J)\hspace{0.2cm}\Max\limits_{(\lambda,U) \in {\cal B}^J} h^J(\lambda,U),$$
where
$${\cal B}^J = \big \{(\lambda, U) \in K^{*0} \times \R^{k\times m} : (L-UA)^T\lambda  \in \R^n_+ \big \}$$ and
$$h^J(\lambda,U) = Ub.$$
This issue was solved by particularizing the general vector Lagrange-type dual introduced in
\cite{JMP}, a vector dual to $(P)$ for which duality statements can be given for every choice of $b\in \R^m$ being obtained, namely
$$(D^L)\hspace{0.2cm}\Max\limits_{(\lambda, z, v) \in {\cal B}^L} h^L(\lambda, z, v),$$
where
$$
{\cal B}^L = \big \{(\lambda,z,v)  \in K^{*0}
 \times  \R^m  \times  \R^k :  \lambda ^Tv - z^Tb\leq 0 \mbox{ and }
 L^T\lambda -A^Tz \in \R^n_+ \big\}$$ and
$$h^L(\lambda,z,v) = v.$$
Recently, in \cite{H5} another vector dual to $(P)$ was proposed for which the duality assertions were shown via complicated set-valued optimization techniques
$$(D^H)\hspace{0.2cm}\Max\limits_{U \in {\cal B}^H} h^H(U),$$
where
$${\cal B}^H = \big \{U \in \R^{k\times m}
:\nexists x\in \R^n_+ \mbox{ such that } (L-UA)x\leq_K 0 \big
\}$$ and
$$h^H(U) = Ub + \Min \big((L-UA)(\R^n_+), K\big).$$

Motivated by the results given in \cite{H5}, we show in this note that a vector dual to $(P)$ already known in the literature (see \cite{WB, BW})
has already ``closed the duality gap in linear vector optimization''. Its weak, strong and converse duality statements hold under the same hypotheses as those for $(D^H)$ and no trouble appears when $b=0$. Moreover, this revisited vector dual to $(P)$ has a vector objective function, unlike $(D^H)$, where the objective function is of set-valued nature, involving additionally solving another vector-minimization problem. So far this vector dual was given only for the case $K=\R^k_+$, but it can be easily extended for an arbitrary nontrivial pointed closed convex cone $K\subseteq \R^k$ as follows
$$
(D)\hspace{0.2cm}\Max\limits_{(\lambda, U, v) \in {\cal
B}} h(\lambda, U, v),$$ where
$$
{\cal B} = \big\{ (\lambda, U, v)  \in K^{*0}\times
\R^{k\times m}\times \R^k: \lambda ^Tv=0 \mbox{ and }
(L-UA)^T\lambda \in \R^n_+\big\} $$and
$$h(\lambda, U, v) = Ub+v.$$

We refer the reader to \cite[Section 5.5]{DVO} for a complete analysis of the relations between all the vector dual problems to $(P)$ introduced above in the case $K=\R^k_+$.

The first new result we deliver in this note is the fact that the efficient and properly efficient solutions to $(P)$ coincide in the considered framework, too, extending the classical statement due to Isermann from \cite{IOR} given for $K=\R^k_+$. Then we prove for $(D)$ weak, strong and converse duality statements. Finally, we show that the relations of inclusions involving the images of the feasible sets through their objective functions of the vector duals to $(P)$ from \cite[Remark 5.5.3]{DVO} remain valid in the current framework, too. Examples when the ``new-old'' vector dual we propose here does not coincide with the other vector duals in discussion are provided, too.

\section{Weak, strong and converse duality for $(P)$ and $(D)$}

An element $\bar x\in {\cal A}$ is said to be a \textit{properly efficient solution (in the sense of linear scalarization)} to $(P)$ if $L\bar x\in \PMin(L({\cal A}), K)$, i.e. there exists $\lambda\in K^{*0}$ such that $\lambda ^T(L\bar x) \leq \lambda ^T(L x)$ for all $x\in {\cal A}$. An element $\bar x\in {\cal A}$ is said to be an \textit{efficient solution} to $(P)$ if $L\bar x\in \Min(L({\cal A}), K)$, i.e. there exists no $x\in {\cal A}$ such that $Lx \leq_K L\bar x$. Of course each properly efficient solution (in the sense of linear scalarization) $\bar x$ to $(P)$ is also efficient to $(P)$. Let us recall now a separation result from the literature, followed by a statement concerning the solutions of $(P)$.\\

{\bf Lemma 1.} (cf. \cite[Lemma 2.2(i)]{H5}) {\it Let $M\subseteq \R^k$ a polyhedral set such that $M\cap K=\{0\}$. Then there exists $\gamma \in \R^k \backslash \{0\}$ such that
$$
\gamma^T k <0 \leq \gamma ^T m\ \forall k\in K\backslash \{0\}\ \forall m\in M.
$$}

{\bf Theorem 1.} {\it Every efficient solution to $(P)$ is properly efficient (in the sense of linear scalarization) to $(P)$.}\\

{\bf Proof.} Let $\bar x\in {\cal A}$ be an efficient solution to $(P)$. Since $\cal A$ is, by construction, a polyhedral set, via \cite[Theorem 19.3]{RC} we get that $L({\cal A})$ is polyhedral, too. Consequently, also $L({\cal A})-L\bar x$ is a polyhedral set. The efficiency of $\bar x$ to $(P)$ yields $(L({\cal A})-L\bar x)\cap (-K) = \{0\}$, thus we are allowed to apply Lemma 1, which yields the existence of $\gamma \in \R^k\backslash \{0\}$ for which
\begin{equation}\label{E1}
\gamma^T (-k) <0 \leq \gamma ^T (Lx-L\bar x)\ \forall k\in K\backslash \{0\}\ \forall x\in {\cal A}.
\end{equation}
Since $\gamma ^Tk>0$ for all $k\in K\backslash \{0\}$, it follows that $\gamma \in K^{*0}$. From \eqref{E1} we obtain $\gamma^T (L\bar x)\leq \gamma^T(Lx)$ for all $x\in {\cal A}$, which, taking into account that $\gamma \in K^{*0}$, means actually that $\bar x$ is a properly efficient solution (in the sense of linear scalarization) to $(P)$.\hfill{$\Box$}$\:$\\

{\it Remark 1.} In Theorem 1 we extend the classical result proven in a quite complicated way in \cite{IOR} for the special case $K=\R^k_+$, showing that the efficient solutions to $(P)$ and its properly efficient solutions (in the sense of linear scalarization) coincide. Actually, this statement remains valid when the feasible set of $(P)$ is replaced by a set $\overline{\cal A}$ for which $L(\overline{\cal A})$ is polyhedral.\\

{\it Remark 2.} In the literature there were proposed several concepts of properly efficient solutions to a vector optimization problem. Taking into account that all these properly efficient solutions are also efficient to the given vector optimization problem and the fact that (see \cite[Proposition 2.4.16]{DVO}) the properly efficient solutions (in the sense of linear scalarization) are properly efficient to
the same problem in every other sense, too, Theorem 1 yields that for $(P)$ the properly efficient solutions (in the sense of linear scalarization)
coincide also with the properly efficient solutions to $(P)$ in the senses of Geoffrion, Hurwicz, Borwein, Benson, Henig and Lampe and generalized Borwein, respectively (see \cite[Section 2.4]{DVO}). Taking into account Theorem 1, it is obvious that it is enough to deal only with the efficient solutions to $(P)$, since they coincide with all the types of properly efficient solutions considered in the literature.\\

Let us show now a Farkas-type result which allows us to formulate the feasible sets of the vector dual problems to $(P)$ in a different manner.\\

{\bf Lemma 2.} {\it Let $U\in \R^{k\times m}$. Then $(L-UA)(\R^n_+)\cap (-K)=\{0\}$ if and only if there exists $\lambda\in K^{*0}$ such that $(L-UA)^T\lambda\in \R^n_+$.}\\

{\bf Proof.}
``$\Rightarrow$'' The set $(L-UA)(\R^n_+)$ is polyhedral and has with the nontrivial pointed closed convex cone $-K$ only the origin as a common element. Applying Lemma 1 we obtain a $\lambda\in \R^k\backslash \{0\}$ for which
\begin{equation}\label{E2}
\lambda ^T (-k) <0 \leq \lambda ^T((L-UA)x)\ \forall x\in \R^n_+\ \forall k\in K\backslash \{0\}.
\end{equation}
Like in the proof of Lemma 1 we obtain that $\lambda\in K^{*0}$ and, by \eqref{E2} it follows immediately that $(L-UA)^T\lambda\in \R^n_+$.

``$\Leftarrow$'' Assuming the existence of an $x\in \R^n_+$ for which $(L-UA)x\in (-K)\backslash \{0\}$, it follows $\lambda^T ((L-UA)x) <0$, but
$\lambda^T ((L-UA)x)= ((L-UA)^T\lambda)^Tx\geq 0$ since $(L-UA)^T\lambda\in \R^n_+$ and $x\in \R^n_+$. The so-obtained contradiction yields $(L-UA)(\R^n_+)\cap (-K)=\{0\}$.\hfill{$\Box$}$\:$\\

Further we prove for the primal-dual pair of vector optimization problems $(P)-(D)$ weak, strong and converse duality statements.\\

{\bf Theorem 2.} (weak duality) {\it There exist no $x\in {\cal A}$ and $(\lambda, U, v) \in {\cal
B}$ such that $Lx\leq_K Ub + v$.}\\

{\bf Proof.} Assume the existence of $x\in {\cal A}$ and $(\lambda, U, v) \in {\cal
B}$ such that $Lx\leq_K Ub + v$. Then $0<\lambda ^T(Ub + v - Lx)= \lambda ^T(U(Ax) - Lx) = -((L-UA)^T\lambda)^Tx \leq 0$, since
$(L-UA)^T\lambda \in \R^n_+$ and $x\in \R^n_+$. But this cannot happen, therefore the assumption we made is false.\hfill{$\Box$}$\:$\\

{\bf Theorem 3.} (strong duality) {\it If $\bar x$ is an efficient solution to $(P)$, there exists $(\bar\lambda, \bar U, \bar v) \in {\cal B}$, an efficient solution to $(D)$, such that $L\bar x= \bar Ub + \bar v$.}\\

{\bf Proof.} The efficiency of $\bar x$ to $(P)$ yields via Theorem 1 that $\bar x$ is also properly efficient to $(P)$. Thus there exists $\bar \lambda\in K^{*0}$ such that $\bar \lambda^T (L\bar x)\leq \bar \lambda^T(Lx)$ for all $x\in {\cal A}$.

On the other hand, one has strong duality for the scalar optimization problem
$$\inf_{x\in {\cal A}}\big\{\bar \lambda^T(Lx)\big\}$$ and its Lagrange dual $$\sup \Big\{-\eta ^Tb: \eta\in \R^m, L^T\bar\lambda + A^T\eta\in \R^n_+\Big\},$$ i.e.
their optimal objective values coincide and the dual has an optimal solution, say $\bar \eta$. Consequently, $\bar \lambda^T (L\bar x) + \bar\eta^Tb=0$ and $L^T\bar\lambda  + A^T\bar \eta\in \R^n_+$.

As $\bar\lambda\in K^{*0}$, there exists $\tilde\lambda \in K\backslash\{0\}$ such that $\bar\lambda^T\tilde\lambda = 1$. Let $\bar U:= - \tilde\lambda\bar \eta^T$ and $\bar v:= L\bar x - \bar Ub$. It is obvious that $\bar U\in \R^{k\times m}$ and $\bar v\in \R^k$. Moreover,
$\bar\lambda^T\bar v= \bar\lambda^T(L\bar x - \bar Ub)= \bar \lambda^T (L\bar x) + \bar\eta^Tb=0$ and $(L-\bar UA)^T\bar \lambda =
L^T\bar\lambda  + A^T\bar \eta\in \R^n_+$. Consequently, $(\bar\lambda, \bar U, \bar v) \in {\cal B}$ and $\bar Ub + \bar v =
\bar Ub + L\bar x - \bar Ub = L\bar x$. Assuming that $(\bar\lambda, \bar U, \bar v)$ is not efficient to $(D)$, i.e. the existence of another feasible solution $(\lambda, U, v) \in {\cal B}$ satisfying $\bar U b+\bar v\leq _K Ub+v$, it follows $L\bar x\leq _K Ub+v$, which contradicts Theorem 2. Consequently, $(\bar\lambda, \bar U, \bar v)$ is an efficient solution to $(D)$ for which  $L\bar x= \bar Ub + \bar v$.\hfill{$\Box$}$\:$\\

{\bf Theorem 4.} (converse duality) {\it If $(\bar\lambda, \bar U, \bar v) \in {\cal B}$ is an efficient solution to $(D)$, there exists
$\bar x\in {\cal A}$, an efficient solution to $(P)$, such that $L\bar x= \bar Ub + \bar v$.}\\

{\bf Proof.} Let $\bar d:=\bar Ub + \bar v$. Assume that ${\cal A}=\emptyset$. Then $b\neq 0$ and, by Farkas' Lemma there exists $\bar z\in \R^m$ such that $b^T\bar z>0$ and $A^T\bar z \in -\R^n_+$. As $\bar\lambda \in K^{*0}$, there exists $\tilde\lambda \in K\backslash \{0\}$ such that $\bar\lambda^T\tilde\lambda = 1$. Let $\tilde U:=\tilde \lambda \bar z^T + \bar U\in \R^{k\times m}$. We have $(L-\tilde U A)^T\bar\lambda = (L-\bar U A)^T\bar\lambda - A^T\bar z\in \R^n_+$, thus $(\bar\lambda, \tilde U, \bar v) \in {\cal B}$. But $h(\bar\lambda, \tilde U, \bar v)= \tilde U b+ \bar v= \tilde\lambda \bar z^Tb+ \bar Ub + \bar v = \tilde\lambda \bar z^Tb + \bar d \geq_K \bar d$, which contradicts the efficiency of $(\bar\lambda, \bar U, \bar v)$ to $(D)$. Consequently, ${\cal A}\neq \emptyset$.

Suppose now that $\bar d\notin L({\cal A})$. Using Theorem 2 it follows easily that $\bar d\notin L({\cal A}) + K$, too.
Since ${\cal A}=A^{-1}(b)\cap \R^n_+$, we have $0^+{\cal A}= 0^+(A^{-1}(b))\cap 0^+\R^n_+$. As $0^+(A^{-1}(b))=A^{-1}(0^+\{b\})= A^{-1}(0)$ and $0^+\R^n_+=\R^n_+$, it follows $0^+{\cal A}= A^{-1}(0)\cap \R^n_+$. Then $0^+ L({\cal A})= L(0^+ {\cal A}) = L(A^{-1}(0)\cap \R^n_+)
= \{Lx: x\in \R^n_+, Ax=0\}\subseteq (L-\bar UA)(\R^n_+)$ and, obviously, $0\in 0^+ L({\cal A})$.

Using Lemma 2 we obtain $(L-\bar UA)(\R^n_+)\cap (-K)=\{0\}$, thus, taking into account the inclusions from above, we obtain $0^+ L({\cal A})\cap (-K)=\{0\} \subseteq K = 0^+K$. This assertion and the fact that $L({\cal A})$ is polyhedral and $K$ is closed convex yield, via \cite[Theorem 20.3]{RC}, that $L({\cal A}) + K$ is a closed convex set. Applying \cite[Corollary 11.4.2]{RC} we obtain a $\gamma\in \R^k\backslash \{0\}$ and an $\alpha\in \R$ such that
\begin{equation}\label{E3}
\gamma ^T \bar d < \alpha < \gamma ^T(Lx + k)\ \forall x\in {\cal A}\ \forall k\in K.
\end{equation}

Assuming that $\gamma\notin K^*$ would yield the existence of some $k\in K$ for which $\gamma^Tk<0$. Taking into account that $K$ is a cone, this implies a contradiction to \eqref{E3}, consequently $\gamma\in K^*$. Taking $k=0$ in \eqref{E3} it follows
\begin{equation}\label{E4}
\gamma ^T \bar d <\alpha < \gamma ^T(Lx)\ \forall x\in {\cal A}.
\end{equation}
On the other hand, for all $x\in {\cal A}$ one has $0\leq \bar \lambda^T((L-\bar UA)x) = \bar \lambda^T(Lx-\bar Ub) =
\bar \lambda^T(Lx-\bar Ub) -  \bar \lambda^T\bar v= \bar \lambda^T (Lx-\bar d)$, therefore
\begin{equation}\label{E5}
\bar \lambda ^T \bar d \leq \bar \lambda ^T(Lx)\ \forall x\in {\cal A}.
\end{equation}
Now, taking $\delta:=\alpha - \gamma^T\bar d>0$ it follows $\bar d^T(s\bar \lambda + (1-s)\gamma)= \alpha - \delta + s(\bar \lambda^T\bar d - \alpha + \delta)$ for all $s\in \R$.
Note that there exists an $\bar s\in (0, 1)$ such that $\bar s(\bar\lambda ^T\bar d - \alpha + \delta)<\delta/2$ and $\bar s(\bar\lambda ^T\bar d - \alpha)>-\delta/2$, and let $\lambda:=\bar s\bar \lambda + (1-\bar s)\gamma$. It is clear that $\lambda\in K^{*0}$.

By \eqref{E4} and \eqref{E5} it follows $s\bar\lambda^T\bar d + (1-s)\alpha < (s\bar \lambda + (1-s) \gamma)^T(Lx)$ for all $x\in {\cal A}$ and all $s\in (0, 1)$, consequently
$$
\lambda ^T \bar d = \bar s\bar \lambda ^T\bar d+ (1-\bar s)\gamma ^T\bar d = \bar s\bar\lambda^T\bar d + (1-\bar s) (\alpha - \delta)
$$
\begin{equation}\label{E6}
<\frac{\delta}{2} + \bar s (\alpha - \delta) + (1-\bar s)(\alpha - \delta)= \alpha-\frac{\delta}{2} < \lambda ^T(Lx)\ \forall x\in {\cal A}.
\end{equation}
Since there is strong duality for the
scalar linear optimization problem $$\inf_{x\in {\cal A}}\big\{\lambda ^T(Lx)\big\}$$ and its Lagrange dual
$$\sup \big\{-\eta ^Tb: \eta\in \R^m, L^T\lambda + A^T\eta\in \R^n_+\big\},$$ the latter
has an optimal solution, say $\bar \eta$, and $\inf_{x\in {\cal A}}\lambda ^T(Lx) +  \bar\eta^Tb=0$ and $L^T\lambda + A^T\bar \eta\in \R^n_+$. As $\bar\lambda\in K^{*0}$, there exists $\tilde\lambda \in K\backslash \{0\}$ such that $\bar\lambda^T\tilde\lambda = 1$. Let $U:= - \tilde\lambda\bar \eta^T$. It follows that $(L-UA)^T\lambda \in \R^n_+$ and $\inf_{x\in {\cal A}}\lambda ^T(Lx) =  \lambda ^T (Ub)$.

Consider now the hyperplane ${\cal H}:= \{Ub + v: \lambda ^Tv=0\}$, which is nothing but the set $\{w\in \R^k: \lambda ^Tw = \lambda ^T(Ub)\}$. Consequently, ${\cal H}\subseteq h({\cal B})$. On the other hand, \eqref{E6} yields $\lambda ^T \bar d < \lambda ^T (Ub)$. Then there exists a $\bar k\in K\backslash \{0\}$ such that $\lambda ^T (\bar d + \bar k) = \lambda ^T (Ub)$, which has as consequence that $\bar d + \bar k\in {\cal H}\subseteq h({\cal B})$.
Noting that $\bar d\leq_K \bar d + \bar k$, we have just arrived to a contradiction to the maximality of $\bar d$ to the set $h({\cal B})$. Therefore
our initial supposition is false, consequently $\bar d\in L({\cal A})$. Then there exists $\bar x\in {\cal A}$ such that $L\bar x=\bar d=\bar Ub + \bar v$. Employing Theorem 2, it follows that $\bar x$ is an efficient solution to $(P)$.\hfill{$\Box$}$\:$\\

{\it Remark 3.} If $\bar x\in {\cal A}$ and $(\bar\lambda, \bar U, \bar v) \in {\cal B}$ are, like in the results
given above, such that $L\bar x= \bar Ub + \bar v$, then the \textit{complementarity condition} $\bar x^T(L-\bar U A)^T\bar\lambda =0$ is fulfilled.\\

Analogously to \cite[Theorem 3.14]{H5} we summarize the results from above in a general duality statement for $(P)$ and $(D)$.\\

{\bf Corollary 1.} {\it One has $\Min(L({\cal A}), K) =  \Max (h({\cal B}), K)$.}\\

To complete the investigation on the primal-dual pair of vector optimization problems $(P)-(D)$ we give also the following assertions.\\

{\bf Theorem 5.} {\it If ${\cal A}\neq\emptyset$, the problem $(P)$ has no efficient solutions if and only if ${\cal B}=\emptyset$.}\\

{\bf Proof.}
``$\Rightarrow$'' By \cite[Lemma 2.1]{H5}, the lack of efficient solutions to $(P)$ yields $0^+L({\cal A}) \cap (-K) \backslash \{0\} \neq\emptyset$. Then $(L-UA)(\R^n_+)\cap (-K)\backslash \{0\} \neq\emptyset$ for all $U\in \R^{k\times m}$ and employing Lemma 2 we see that ${\cal B}$ cannot contain in this situation any element.

``$\Leftarrow$'' Assuming that $(P)$ has efficient solutions, Theorem 3 yields that also $(D)$ has an efficient solution. But this cannot happen since the dual has no feasible elements, consequently $(P)$ has no efficient solutions.\hfill{$\Box$}$\:$\\

{\bf Theorem 6.} {\it If ${\cal B}\neq\emptyset$, the problem $(D)$ has no efficient solutions if and only if ${\cal A}=\emptyset$.}\\

{\bf Proof.}
``$\Rightarrow$'' Assume that ${\cal A}\neq\emptyset$. If $(P)$ has no efficient solutions, Theorem 5 would yield ${\cal B}=\emptyset$, but this is false, therefore $(P)$ must have at least an efficient solution. Employing Theorem 3 it follows that $(D)$ has an efficient solution, too, contradicting the assumption we made. Therefore ${\cal A}=\emptyset$.

``$\Leftarrow$'' Assuming that $(D)$ has an efficient solution, Theorem 4 yields that $(P)$ has an efficient solution, too. But this cannot happen since this problem has no feasible elements, consequently $(D)$ has no efficient solutions.\hfill{$\Box$}$\:$\\

\section{Comparisons to other vector duals to $(P)$}

As we have already mentioned in the first section (see also \cite[Section 4.5 and Section 5.5]{DVO}), several vector dual problems to $(P)$ were proposed in the literature and for them weak, strong and converse duality statements are valid under different hypotheses.
For $(D^I)$ (in case $K=\R^k_+$) and $(D^J)$ weak duality holds in general, but for strong and converse duality one needs to impose additionally the condition $b\neq 0$ to the hypotheses of the corresponding theorems given for $(D)$. For $(D^H)$ all three duality statements hold under the
hypotheses of the corresponding theorems regarding $(D)$. Concerning $(D^L)$, weak and strong duality were shown (for instance in \cite{DVO}), but the converse duality statement does not follow directly, so we prove it.\\

{\bf Theorem 7.} {\it If $(\bar\lambda, \bar z, \bar v) \in {\cal B}^L$ is an efficient solution to $(D^L)$, there exists
$\bar x\in {\cal A}$, an efficient solution to $(P)$, such that $L\bar x= \bar v$.}\\

{\bf Proof.} Analogously to the proof of Theorem 4 it can be easily shown that ${\cal A}\neq\emptyset$. Since $\bar\lambda\in K^{*0}$, there exists $\tilde\lambda \in K\backslash\{0\}$ such that $\bar\lambda^T\tilde\lambda = 1$.
Let $U:= (\bar z\tilde\lambda^T)^T$. Then $U^T\bar\lambda = \bar z\tilde\lambda^T\bar\lambda = \bar z$. Thus, $(L-UA)^T\bar\lambda = L^T\bar\lambda - A^TU^T\bar\lambda = L^T\bar\lambda - A^T\bar z\in \R^n_+$. Assuming the existence of some $x\in \R^n_+$ for which $(L-UA)x\in - K\backslash \{0\}$, it follows $\bar \lambda ^T ((L-UA)x) <0$. But $\bar \lambda ^T ((L-UA)x)= x^T ((L-UA)^T\bar\lambda) \geq 0$, since $x\in \R^n_+$ and $(L-UA)^T\bar\lambda\in \R^n_+$. This contradiction yields $(L-UA)(\R^n_+)\cap (-K) = \{0\}$. Like in the proof of Theorem 4, this result, together with the facts that $L({\cal A})$ is polyhedral and $K$ is closed convex implies, via \cite[Theorem 20.3]{RC}, that $L({\cal A}) + K$ is a closed convex set. The existence of $\bar x\in {\cal A}$ properly efficient, thus also efficient, solution to $(P)$ fulfilling $L\bar x= \bar v$ follows
in the lines of \cite[Theorem 4.3.4]{DVO} (see also \cite[Section 4.5.1]{DVO}). \hfill{$\Box$}$\:$\\

Let us see now what inclusions involving the images of the feasible sets through their objective functions of the vector duals to $(P)$
considered in this paper can be established. In \cite{H5} it is mentioned that $h^J({\cal B}^J) \subseteq h^H({\cal B}^H)$, an example when
these sets do not coincide being also provided. For $(D)$ and $(D^H)$ we have the following assertion.\\

{\bf Theorem 8.} {\it It holds $h^H({\cal B}^H) \subseteq h({\cal B})$.}\\

{\bf Proof.} Let $d\in h^H({\cal B}^H)$. Thus, there exist $\bar U\in {\cal B}^H$ and an efficient solution $\bar x\in \R^n_+$ to
\begin{equation}\label{P1}
\Min_{x\in \R^n_+} \{(L-\bar UA)x\},
\end{equation} such that $d=h^H(\bar U)=\bar Ub + (L-\bar UA)\bar x$.

The efficiency of $\bar x$ to the problem \eqref{P1} yields, via Theorem 1, that $\bar x$ is a properly efficient solution to this problem, too. Consequently, there exists $\gamma \in K^{*0}$ such that
\begin{equation}\label{E7}
\gamma^T((L-\bar UA)\bar x) \leq \gamma^T((L-\bar UA) x)\ \forall x\in \R^n_+.
\end{equation}
This yields $\gamma^T((L-\bar UA)\bar x) \leq 0$. On the other hand, taking in \eqref{E7} $x:=x+\bar x\in \R^n_+$ it follows immediately
$\gamma^T((L-\bar UA) x)\geq 0$ for all $x\in \R^n_+$. Therefore $\gamma^T((L-\bar UA)\bar x) \geq 0$, consequently $\gamma^T((L-\bar UA)\bar x) = 0$. Taking $\bar v= (L-\bar UA)\bar x$, it follows $\gamma^T\bar v=0$ and, since $\gamma^T(L-\bar U A)\in \R^n_+$, also $(\gamma, \bar U, \bar v)\in {\cal B}$. As $d=h^H(\bar U)=\bar Ub + (L-\bar UA)\bar x = \bar Ub + \bar v = h(\gamma, \bar U, \bar v)\in h({\cal B})$, we obtain
 $h^H({\cal B}^H) \subseteq h({\cal B})$.\hfill{$\Box$}$\:$\\

{\it Remark 4.} An example showing that the just proven inclusion can be sometimes strict was given in \cite[Example 5.5.1]{DVO}.\\

{\bf Theorem 9.} {\it It holds $h({\cal B}) \subseteq h^L({\cal B}^L)$.}\\

{\bf Proof.} Let $d\in h({\cal B})$. Thus, there exist $(\lambda, U, v)\in {\cal B}$ such that $d=h(\lambda,  U,  v)= Ub +  v$.

Let $ z:= U^T\lambda$. Then $\lambda ^Td = \lambda^T(Ub + v) = (\lambda^T(Ub + v))^T = b^T(U^T\lambda) + v^T\lambda = b^Tz$, while
$L^T\lambda- A^Tz = L^T\lambda - A^T(U^T\lambda)= (L-UA)^T\lambda \in \R^n_+$. Consequently, $(\lambda, z, d)\in {\cal B}^L$ and, since
$d= h(\lambda, U, v) = h^L(\lambda, z, d)\in h^L({\cal B}^L)$, it follows that $h({\cal B}) \subseteq h^L({\cal B}^L)$.\hfill{$\Box$}$\:$\\

{\it Remark 5.} To show that the inclusion proven above does not turn into equality in general, consider the following situation. Let $n=1$, $k=2$, $m=2$, $L=(0, 0)^T$, $A=(1, 1)^T$, $b=(-1, -1)^T$ and $K=\R^2_+$. Then, for $v=(-1, -1)^T$, $\lambda=(1, 1)^T$ and $z=(0, 0)^T$ we have $(\lambda, z, v)\in {\cal B}^L$ since $\lambda^Tv=-2\leq 0 = z^Tb$ and $L^T\lambda - A^Tz = (0, 0)^T\in \R^2_+$. Consequently, $h^L(\lambda, z, v) = (-1, -1)^T\in h^L({\cal B}^L)$.

On the other hand, assuming that $(-1, -1)\in h({\cal B})$, there must exist $(\bar \lambda, \bar U, \bar v)$ $\in {\cal B}$ such that
$h(\bar \lambda, \bar U, \bar v)= \bar Ub + \bar v= (-1, -1)^T$. Then $\bar \lambda ^T(\bar Ub+\bar v)= - \bar \lambda_1 - \bar \lambda_2 <0$, where
$\bar \lambda=(\bar \lambda_1, \bar \lambda_2)^T\in (\R^2_+)^{*0}=\inte (\R^2_+)$. But $\bar \lambda ^T(\bar Ub+\bar v)= \bar \lambda ^T\bar U (-1, -1)^T= - (\bar U A)^T\bar \lambda = (L-\bar UA)^T\bar \lambda \geq 0$, which contradicts what we obtained above as a consequence of the assumption we made. Therefore $(-1, -1)\notin  h({\cal B})$.

\section{Conclusions}

Taking into consideration what we have proven in this note, one can conclude that the images of the feasible sets through their objective functions of the vector duals to $(P)$ we dealt with respect the following inclusions chain
$$h^J({\cal B}^J) \subsetneqq h^H({\cal B}^H) \subsetneqq h({\cal B}) \subsetneqq h^L({\cal B}^L),$$
while the the sets of maximal elements of these sets fulfill
$$\Max(h^J({\cal B}^J), K) \subsetneqq \Min(L({\cal A}), K) = \Max (h^H({\cal B}^H), K)$$
$$= \Max (h({\cal B}), K) = \Max (h^L({\cal B}^L), K).$$
Thus the schemes from \cite[Remark 5.5.3]{DVO} remain valid when the vector-minimization/ maximization is considered with respect to a nontrivial pointed closed convex cone $K\subseteq \R^k$.

Therefore, we have brought into attention with this note a valuable but neglected vector dual to the classical linear vector optimization problem. The image of the feasible set through its objective function of this ``new-old'' vector dual lies between the ones of the vector duals introduced in \cite{H5} and \cite{JMP}, respectively, all these three sets sharing the same maximal elements. Different to the vector dual introduced in \cite{H5}, the dual we consider is defined without resorting to complicated constructions belonging to set-valued optimization and the duality statements regarding it follow more directly. Solving our vector dual problem does not involve the determination of the sets of efficient solutions of two vector optimization problems, as is the case for the vector dual from \cite{H5}. Different to the Lagrange-type vector dual from \cite{JMP}, the revisited vector dual we consider preserves the way the classical vector dual problem due to Isermann and the dual abstract linear optimization problem were formulated, improving their duality properties by not failing like them when $b=0$.

Moreover, we have extended a classical result due to Isermann, showing that the efficient solutions to
the classical linear vector optimization problem coincide with the properly efficient elements in any sense to it also when
the vector-minimization is considered with respect to a nontrivial pointed closed convex cone $K\subseteq \R^k$ instead of $\R^k_+$.

 \end{document}